\documentclass[11pt]{article}

\usepackage{amscd,amsmath,amsthm,latexsym}
\usepackage{amsfonts}
\usepackage{amssymb}
\usepackage{longtable}
\usepackage[all]{xy}
\usepackage{graphicx}

\newenvironment{pf}{\medskip\noindent{\it Proof.}\enspace}
{\hfill \newline \smallskip}

\numberwithin{equation}{section}

\def\eps{\varepsilon}
\newcommand{\af}{\alpha}
\newcommand{\bt}{\beta}

\newcommand{\dt}{\delta}
\newcommand{\ep}{\varepsilon}

\newcommand{\et}{\eta}

\newcommand{\ld}{\lambda}

\newcommand{\kp}{\kappa}

\newcommand{\ift}{\infty}
\newcommand{\ds}{\displaystyle}

\def\1{{\bf 1}}

\def\sA {{\cal A}}  \def\sC {{\cal C}}

  \def\sL {{\cal L}}

\def\R {{\mathbb R}} \def\RR {{\mathbb R}}

\newcommand{\expr}[1]{\left( #1 \right)}

\newtheorem{thm}{Theorem}[section]
\newtheorem{lemma}[thm]{Lemma}
\newtheorem{defn}[thm]{Definition}

\newtheorem{prop}[thm]{Proposition}

\newtheorem{cor}[thm]{Corollary}

\numberwithin{equation}{section}

\def\qed{{\hfill $\Box$ \bigskip}}

\def\RR{{\mathbb R}}
\def\R{{\mathbb R}}

\def\S{{\bf S}}
\def\A{{\bf A}}
\def\E{{\mathbb E}}

\def\P{{\mathbb P}}

\def\eps{\varepsilon}

\def\pf{\noindent{\bf Proof.} }
\def\beq{\begin{equation}}
\def\eeq{\end{equation}}
\def\bee{\begin{equation}}
\def\eee{\end{equation}}

\begin{document}

\bibliographystyle{plain}

\title{\bf Oscillation of harmonic functions for subordinate Brownian motion and its applications
}
\author{{\bf Panki Kim}\thanks{This research was supported by Basic Science Research Program through the National Research Foundation of Korea(NRF) funded by the Ministry of Education, Science and Technology (0409-20120034)
} \quad and  \quad {\bf Yunju Lee}}

\date{(September 27, 2012)}
\maketitle

\begin{abstract}
In this paper, we establish an oscillation estimate of nonnegative harmonic functions for a pure-jump subordinate Brownian motion. The infinitesimal generator of such subordinate Brownian motion is an integro-differential operator. As an application, we give a probabilistic proof of the following form of relative Fatou theorem for such subordinate Brownian motion $X$ in bounded $\kappa$-fat open set; if $u$ is a positive harmonic function with respect to $X$ in  a bounded $\kappa$-fat open set $D$ and $h$ is a positive  harmonic function in $D$ vanishing on $D^c$, then the non-tangential limit
 of  $u/h$ exists almost everywhere with respect to the Martin-representing measure of $h$.
\end{abstract}
\vspace{.6truein}



\section{Introduction}\label{sec-intro}
Nowadays L\'evy processes
have been receiving intensive study due to their importance both in theories and applications.
They
are widely used in various fields, such as
mathematical finance, actuarial mathematics
and  mathematical physics. Typically, the infinitesimal generators of general L\'evy processes
 in $\R^d$ are not differential
operators but  integro-differential operators. Even though  integro-differential operators are very
important in the theory of partial differential equations, general L\'evy processes and corresponding integro-differential operators  are not  easy to deal with.
 For a summary of some of these
recent results from the probability literature, one can see
\cite{BBKRSV} and the references therein. We refer readers
to \cite{CSS, CSi} for samples of
 recent progresses in the PDE literature.

Let $W=(W_t:\, t\ge 0)$ be a Brownian motion in
$\R^d$ and $S=(S_t:\, t\ge 0)$ be a subordinator independent of $W$.
The process $X=(X_t:\, t\ge 0)$ defined by
$X_t=W_{S_t}$ is a rotationally invariant L\'evy process in $\R^d$ and is called a
subordinate Brownian motion.
Subordinate Brownian motions form a very large class of L\'evy processes. Nonetheless, compared with general L\'evy processes, subordinate
Brownian motions are much more tractable. If we take the Brownian
motion
$W$ as given, then $X$ is completely determined by the Laplace exponent of
subordinator $S$. Hence one can deduce the properties of $X$ from the subordinator $S$, or equivalently the Laplace exponent of it.

The purpose of
 this paper is to give an oscillation
 estimate for (unbounded) harmonic functions (see Section \ref{sec-prelim} for the definition of harmonicity)  for a large class of subordinate Brownian motions. Then using our estimates, we discuss non-tangential limits of the ratio of two harmonic functions with respect to such subordinate Brownian motions.

Now we state the first main result of this paper.

\begin{thm} \label{th:osc}
Suppose that $X=(X_t:\, t\ge 0)$ is a L\'evy process  whose
characteristic exponent is given by $\Phi(\theta)=\phi(|\theta|^2)$,
$\theta\in \R^d$, where
$
\phi:(0,\infty)\to [0,\infty)$  is a complete Bernstein
function such that   $ \phi(\lambda)=
\lambda^{\alpha/2}\ell(\lambda)$, $\alpha\in (0, 2)$ and $\ell:(0,\infty)\to (0,\infty)$ is slowly varying  at $\ift$.
Then for every $\eta > 0$, there exists $a=a(\eta, \af, d, \ell) \in (0, 1)$  such that for every $x_0 \in \R^d$ and $r \in (0, 1]$,
\begin{equation*}
  \sup_{x\in B(x_0,ar)} u(x) \leq
  (1 + \eta){\inf_{x\in B(x_0,ar)} u(x)}
\end{equation*}
for every nonnegative function $u$ in $\RR^d$
which is harmonic in $B(x_0, r)$ with respect to $X$. \end{thm}
Note that, for unlike a local operator, Theorem \ref{th:osc} can not be obtained from  Harnack inequality and Moser's iteration method because harmonic functions in
Theorem \ref{th:osc}
are
nonnegative  in  the whole space $\RR^d$.
On the other hand, if one just assumes that a harmonic function is
nonnegative  in $B(x_0, 2r)$, then even Harnack inequality does not hold (see \cite{Ka2}).

Recently many results are obtained under the weaker assumption that $\phi$ is comparable to a regularly varying function at $\ift$ (see \cite{KPS, KSV2, KSV3, KSV4}). But our technical Lemmas \ref{l:j1}--\ref{l:p1} cannot be obtained under such assumptions.

Doob  proved  the relative Fatou  theorem in the classical sense (\cite{D3}). That is,
the ratio $u/h$ of two positive harmonic functions with respect to Brownian motion on a unit open ball has non-tangential limits almost everywhere with respect to the Martin measure of $h$. Later, relative Fatou  theorem in the classical sense has been extended to some general open sets (see \cite{W} and references therein). But relative Fatou  theorem stated above and Fatou  theorem are not true for harmonic functions for the fractional Laplacian $\Delta^{\alpha/2} : = - (-\Delta^{\alpha/2}) $ when $\alpha \in (0,2)$ (see \cite{BY} for some counterexamples). Correct formulation of relative Fatou  theorem for integro-differential operator is
 the existence of non-tangential limits of  the ratio $u/h$, where $u$ is positive harmonic in a open set $D$ and $h$ is a positive
harmonic function in $D$ vanishing on $D^c$ (see \cite{BD, K2, KS, MR}).

 In this paper, through a probabilistic method and Theorem \ref{th:osc}, we show
 in Theorem \ref{T:Fatou}
 that relative Fatou  theorem holds for subordinate Brownian motion in very general open sets, namely, bounded $\kp$-fat open sets, the family that includes bounded Lipschitz open sets.

This paper is organized as follows. In section 2, we recall the definition of subordinate Brownian motion and its basic properties
under our assumptions.
In Section 3, we give the proof of Theorem \ref{th:osc}. In these sections, the influence of \cite{BKK} in our results will be apparent.
Section 4 contains the proof of relative Fatou  theorem
in  bounded $\kappa$-fat open sets.
The main idea of our proof is similar to \cite{K2}, which is inspired by Doob's approach (see also \cite{B}). We use Harnack and boundary Harnack principle obtained in \cite{KSV} and our Theorem \ref{th:osc}. If the open set is the unit ball in $\RR^2$, we show that our result is the best possible one.

In the sequel, we will use the following convention: The
value of the constant $C_*$  will remain
the same throughout this paper, while the constants
$c_0, c_1, c_2, \cdots$ signify  constants whose values are
unimportant and which may change from location to location. The labeling of the
constants $c_0, c_1, c_2, \cdots$ starts anew in the statement of each
result. We use ``$:=$" to denote a definition, which
is  read as ``is defined to be". We denote $a \wedge b := \min \{ a, b\}$, $a \vee b := \max \{ a, b\}$ and $f(t) \sim g(t)$, $t \to 0$ ($f(t) \sim g(t)$, $t \to \infty$, respectively) means $
\lim_{t \to 0} f(t)/g(t) = 1$ ($\lim_{t \to \infty} f(t)/g(t) = 1$, respectively).
   For any open set $U$, we denote $\delta_U(x)=\text{dist} (x, U^c)$. Let  $
A(x, a, b) :=\{ y \in \R^d: a \le  |x-y| <b\}$ and $B(x_0, r)$ be a ball in $\R ^d$ centered at $x_0$ whose radius is $r$.
When $x_0$ is the origin, we simply denote $B_r:=B(0,r)$.

\section{Preliminaries} \label{sec-prelim}

Suppose that $S=(S_t: t\ge 0)$ is a subordinator, that is, an increasing
L\'evy process taking values in $[0,\infty)$ with $S_0=0$. A
subordinator $S$ is completely characterized by its Laplace exponent
$\phi$ via
$$
\E[\exp(-\lambda S_t)]=\exp(-t \phi(\lambda))\, ,\quad  \lambda > 0.
$$

 A smooth function $\phi : (0, \infty) \to [0, \infty)$ is called a Bernstein function if $(-1)^n D^n \phi \le 0$ for every natural number $n$.
 Every Bernstein function has a representation
$$
\phi (\lambda)=a+b\lambda +\int_{(0,\infty)}(1-e^{-\lambda t})\,
\mu(dt)
$$
where $a,b\ge 0$ and $\mu$ is a measure on $(0,\infty)$ satisfying
$\int_{(0,\infty)}(1\wedge t)\, \mu(dt)<\infty$. $a$ is called the
killing coefficient, $b$ is the drift and $\mu$ is the L\'evy measure of
the Bernstein function. A nonnegative function $\phi$ on $(0,
\infty)$ is the Laplace exponent of a subordinator if and only if it
is a Bernstein function with $\phi(0+)=0$. We also call $\mu$ the L\'evy measure
of the subordinator $S$. A Bernstein function $\phi$ is called a complete Bernstein function
if $\mu$ has a completely monotone density
$t \mapsto \mu(t)$, i.e., $\mu (t) dt = \mu (dt)$ and $(-1)^n D^n \mu\ge 0$ for every non-negative integer
$n$.

Throughout this paper we will assume that\\
\noindent
{\bf (A1) } : {\it
$\phi$ is a complete Bernstein function and regularly varying of index $\alpha/2$ at $\infty$ for some $\alpha\in (0, 2)$. That is, 
\begin{equation}\label{e:asumptiononle}
\phi(\lambda)=\lambda^{\alpha/2}\ell(\lambda)
\end{equation}
for some $\alpha\in (0, 2)$ and some positive function $\ell$ which
is slowly varying at $\infty$.}\\
Note that, this is an assumption about $\phi$ at $\ift$ and nothing is assumed about the behavior near zero. 
Clearly \eqref{e:asumptiononle} implies that $b=0$ and
 $\lambda \to \ell(\lambda)$ is strictly positive and continuous on $(0, \infty)$. We refer to \cite{KSV} for examples. From \cite[Proposition 5.23]{BBKRSV}, we get
\begin{equation}\label{e:abofldofSat0}
\mu(t)\sim  \frac{\alpha}{2\Gamma(1-\alpha/2)} t^{-1}\phi(t^{-1})
 \qquad \text{as } \,\,\, t\to 0
\end{equation}
where $\Gamma (\ld) := \int_0^{\ift} t^{\ld -1} e^{-t} dt $.

 Let $W:=(W_t,\, \P_x  :  t\ge 0, \, x \in \R^d)$ be a Brownian motion on $\R^d$ with $\P_x (W_0 = x) =1 $ and $ \ds
\E_x [e^{i\xi\cdot(W_t-W_0)} ]=
e^{-t|\xi|^2} \, \text{ for } \xi \in \R^d, \,t>0 \text{ and } x \in \R^d.
$
In the remainder of this paper we will use $X=(X_t, \, \P_x  :  t\ge 0, \, x \in \R^d)$ to denote
the subordinate Brownian motion defined by $X_t=W_{S_t}$, where $S= (S_t, t \ge 0)$ is a subordinator whose Laplace exponent is $\phi$ and $S$ is independent of $W$.

Let
\begin{equation}\label{j11}
j(r):=\int^{\infty}_0(4\pi t)^{-d/2}e^{-r^2/(4t)}\mu(t)dt \qquad\text{for }\,\, r>0
\end{equation}
where $\mu(t)$ is the L\'evy density of $S$. Then $J(x) := j(|x|)$ is the L\'evy density of $X$.
Note that the function $r\mapsto
j(r)$ is strictly positive, continuous and decreasing on $(0, \infty)$.
Since $| {\partial}/{\partial r} (e^{-r^2/(4t)} )| \,=\, 4  r^{-1} \big({r^2}/{(8t)}\, e^{-r^2/(8t)}\big) e^{-r^2/(8t)}  \,\le\,  c\,r^{-1}e^{-r^2/(8t)}$ and $ \int^{\infty}_0(4\pi t)^{-d/2}  r^{-1}e^{-r^2/(8t)}\mu(t)dt  \,=\,r^{-1}\,j({r}/{\sqrt2})$, $j'(r)$ is well-defined and is continuous.

Applying \cite[Lemma 13.3.1]{KSV2}, we have the following.

\begin{thm}\label{t:Gorigin}
$$
j(r)\,\sim\, \frac{\alpha\Gamma((d+\alpha)/2)}{2^{1-\alpha}\pi^{d/2}\Gamma(1-
\alpha/2)}\frac{\phi(r^{-2})}{r^{d}}
\qquad\text{as }\,\, r\to 0.
$$
\end{thm}

As an immediate consequence of Theorem \ref{t:Gorigin} and the continuity of $r \mapsto j(r)$ on $(0, \infty)$, we have

\begin{cor}\label{l:J} For every $R>0$, there exists $c=c(R, \af,d, \ell) >1$ such that for every positive
$y$ with $|y|\le R$,
$$
c^{-1} |y|^{-d} \phi(|y|^{-2}) \,\le\, J(y) \,\le\, c \,|y|^{-d} \phi(|y|^{-2}).
$$

\end{cor}

By \cite[Proposition 13.3.5]{KSV2}, the function $r \mapsto j(r)$ enjoys the following properties.
\begin{prop}
\begin{description}
\item{(1)} For any $M>0$, there exists $c_1=c_1(M)>0$ such that
$
j(r)\,\le\, c_1 j(2r)$ for every $r\in (0, M)$.
\item{(2)} There exists $c_2>0$ such that
$
j(r)\,\le\, c_2 j(r+1)$ for every $ r>1$.
\end{description}
\end{prop}
For any open set $D$, we use $\tau_D$ to denote the first exit
time of $D$, i.e., $\tau_D=\inf\{t>0: \, X_t\notin D\}$.
Given  an open set $D\subset \R^d$, we define
$X^D_t(\omega)=X_t(\omega)$ if $t< \tau_D(\omega)$ and
$X^D_t(\omega)=\partial$ if $t\geq  \tau_D(\omega)$, where
$\partial$ is a cemetery state. We now recall the definition
of harmonic functions with respect to $X$.

\begin{defn}\label{def:har1}
Let $D$ be an open subset in $\R^d$.
A function $u$ defined on $\R^d$ is said to be

\begin{description}
\item{\rm{(1)}}  harmonic in $D$ with respect to $X$ if
$
\E_x\left[|u(X_{\tau_{B}})|\right] <\infty
$  and $
u(x)= \E_x[u(X_{\tau_{B}})]$ for every $x\in B$ and open set $B$ whose closure is a compact
subset of $D$;

\item{\rm{(2)}}
regular harmonic in $D$ with respect to $X$ if it is harmonic in $D$
with respect to $X$ and
for each $x \in D$,
$
u(x)= \E_x\left[u(X_{\tau_{D}})\right];
$

\item{\rm{(3)}} harmonic with respect to $X^D$ if it is harmonic with respect to $X$ in $D$ and
vanishes outside $D$.

\end{description}
\end{defn}

By \cite[Corollary 13.4.8]{KSV2}, we have
the following Harnack inequality.

\begin{thm}\label{HP} \emph{(Harnack inequality)}
There exists a constant $C_0>0$  such that for  every $r \in (0,1), \,\, x_0 \in \R^d$ and
function $f \ge 0$ in $\R^d$ which is harmonic in $B(x_0, r)$ with respect to $X$,  we
have
$$
\sup_{y \in B(x_0, r/2)}f(y) \le C_0 \inf _{y \in B(x_0, r/2)}f(y).
$$
\end{thm}

It follows from \cite[Chapter 5]{BBKRSV} that the process $X$ has a transition
density $p(t, x, y)$ which is jointly continuous. By the joint continuity and the strong Markov property, one can easily check that
$$
p_D(t, x, y)\,:=\,p(t, x, y)\,-\,\E_x[ \,p(t-\tau_D, X_{\tau_D}, y)\,;\, t>\tau_D] \qquad \text{for }\,x, y \in D
$$
is the transition density of $X^D$, which is jointly continuous (for example, see \cite[Lemma 5.5]{KPS}).
 For any bounded open set $D\subset \R^d$, we
will use $G_D$ to denote the Green function of $X^D$, i.e.,
$$
G_D(x, y):=\int^\infty_0p_D(t, x, y)dt \qquad \text{for }\,\,x, y\in D.
$$
Note that $G_D$ is continuous in $(D\times D)\setminus\{(x, x): x\in D\}$.

We define the Poisson kernel $P_D (x,y)$ as \begin{equation*}
P_D(x,y)\,:=  \int_{D}
G_D(x,z) J(z-y) \,dz \,  \qquad \text{for }\,(x,y) \in \R^d \times
  {\overline{D}}^c.
\end{equation*}
Thus we have for every
bounded open subset $D$, function
$f \ge 0$ and $x \in D$,
\begin{equation}\label{newls}
\E_x\left[f(X_{\tau_D});\,X_{\tau_D-} \not= X_{\tau_D}  \right]
=\int_{\overline{D}^c} P_D(x,y)f(y)dy.
\end{equation}
Using the continuities of $G_D$ and $J$, one can easily check
that $P_D$ is continuous on $D \times
\overline{D}^c$. Moreover, from \cite[Theorem 1]{Sz1} we know $\P_x (X_{\tau_{B_r}} \in \partial B_r) = 0  $ for $x \in B_r$. Thus every harmonic function $u$ in $D$ is written as 
\bee\label{e:harmonic}
u(x) = \int_{B_r^{c}} P_{B_r} (x,y) u(y) dy \, \qquad \text{for }x\in B_r \subset \overline{B_r} \subset D \,.
\eee
When $r  \le 1$, by the continuity of $P_{B(x_0, r)}$ and Harnack inequality (Theorem \ref{HP}), we get 
$$P_{B(x_0, r)} (x,y) \,\le\, C_0 \, P_{B(x_0, r)} (x_0, y) \quad  \text{ for every } (x,y) \in B(x_0, r/2)\times \overline{B(x_0,r)} ^{\,c}.$$ Since $P_{B(x_0, r)} (x_0 ,y) |u(y)| \in L^1(D)$ for $y \in \overline{B(x_0,r)} ^{\,c}$ by the definition of the harmonicity, applying   Lebesgue dominated convergence theorem to \eqref{e:harmonic} we see that every harmonic function in $D$ with respect to $X$ is continuous.

\section{Oscillation of harmonic functions}

Recall that $S_t$ is a subordinator with Laplace exponent $\phi$, $W$ is a Brownian motion independent with $S_t$ and $X_t = W_{S_t}$. First we show that $\phi$ being a complete Bernstein function
implies that its L\'evy density of $X$ cannot decrease too fast in the
following sense:

\begin{lemma}\label{l:mu}
\begin{equation*}
\lim_{\delta \downarrow 0} \sup_{t >1}\frac{\mu(t)}{\mu(t+\delta)} = 1.
\end{equation*}
\end{lemma}

\pf
Let $\eta>0$ be given.
Since $\mu$ is a completely monotone function, by Bernstein's
theorem (\cite[Theorem 1.4]{SSV}) there exists a measure $m$ on
$[0,\infty)$ such that $\mu(t)=\int_{[0,\infty)}e^{-tx} m(dx).$
Choose $r=r(\eta)>0$ such that
$$ \eta \int_{[0, r]}e^{-x}\, m(dx)\ge \int_{(r,
\infty)}e^{-x}\, m(dx).$$ Then for any $t>1$, we have
\begin{eqnarray*}
&&\eta\int_{[0, r]}e^{-t x}\, m(dx)=\eta\int_{[0, r]}e^{-(t-1) x}e^{-x}\, m(dx)\ge e^{-(t -1)r} \eta \int_{[0, r]}e^{-x}\, m(dx)\\
&&\ge e^{-(t -1)r}\int_{(r, \infty)}e^{-x}\,  m(dx)=\int_{(r, \infty)}e^{-(t -1)r} e^{-x}\,  m(dx)\,\ge \,
\int_{(r, \infty)}e^{-t x}\, m(dx).
\end{eqnarray*}
Thus for any $t>1$ and $\delta>0$,
\begin{eqnarray*}
&&\mu(t+\delta)\ge \int_{[0, r]}e^{-(t+\delta) x}\, m(dx)\ge  e^{-r\delta}\int_{[0,
r]}e^{- t x}\, m(dx)\\
&&=e^{-r\delta} (1+ \eta)^{-1} \left(\int_{[0,
r]}e^{- t x}\, m(dx) +\eta \int_{[0,
r]}e^{- t x}\, m(dx) \right)\\
&& \ge e^{-r\delta} (1+ \eta)^{-1} \left(\int_{[0,
r]}e^{- t x}\, m(dx) +\int_{(r, \infty)}e^{-t x}\, m(dx) \right)\\
&&=      e^{-r\delta}(1+ \eta)^{-1} \int_{[0, \infty)}e^{-t
x}\, m(dx)=\,  e^{-r\delta}(1+ \eta)^{-1}\mu(t).
\end{eqnarray*}
Therefore,
$$
\limsup_{\delta \downarrow 0} \left( \sup_{t >1}\frac{\mu(t)}{\mu(t+\delta)} \right) \,\le\, 1+ \eta\, .
$$
Since $\eta>0$ is arbitrary and $\frac{\mu(t)}{\mu(t+\delta)} \ge 1$, we conclude that this lemma holds.
\qed

\begin{lemma}\label{l:j1}
\begin{equation*}
\lim_{\delta \downarrow 0} \sup_{r >2}\frac{j(r)}{j(r+\delta)} = 1.
\end{equation*}
\end{lemma}

\pf
Fix $\eps \in (0,1)$ and let $L:= \frac{\alpha}{2\Gamma(1-\alpha/2)} $. Using \eqref{e:asumptiononle}, \eqref{e:abofldofSat0}  and the fact that $\ell$ is slowly varying, we choose $t_* = t_* (\eps) \in (0,1/2)$ such that  for every $t \le 2\, t_*$,
\bee \label{e:mu01}
(1+\eps)^{-1}
L \frac{\phi(t^{-1})}{
t} \, \le \, \mu(t) \le (1+\eps)\,
L\frac{\phi(t^{-1})}{
t}
\quad \text{ and } \quad
  1 \,\le \,
  \frac{\phi\big((1+\eps)t^{-1}\big)}
  {\phi(t^{-1})} \,\le\, (1+\eps)^{1+\af/2} .
\eee
By \eqref{e:mu01}
we get
\begin{eqnarray}
\mu\big((1+\ep)t\big) &\ge& (1+\ep)^{-1} \, L \,\frac{\phi\big((1+\ep)^{-1} t^{-1}\big)}{(1+\ep) t } \,\ge\,  (1+\ep)^{-3-\af/2} \,L\, \frac{\phi(t^{-1})}{t} \nonumber \\
&\ge&  (1+\ep)^{-4-\af/2}\, \mu(t) \quad \text{ for every  } t \le 2\,t_* .  \label{e:mu02}
\end{eqnarray}
Now using Lemma \ref{l:mu}, we choose $\delta_1 \in (0, \eps (1+\eps)^{-1}]$ such that for every $t \ge 1$,
\bee \label{e:mu03}
 \mu(t+\delta_1) \le \mu(t)  \le (1+\eps) \mu(t+\delta_1) .
 \eee
Since $$\frac{\mu(t)-\mu\big((1-\delta)^{-1}t\big)}{\mu\big((1-\delta)^{-1}t\big)} \,\le\, \frac{\mu(t)-\mu\big((1-\delta)^{-1}t\big)}{\mu(4)} \quad \text{and} \quad \frac{\mu(t)-\mu(\delta+t)}{\mu(\delta+t)} \,\le\, \frac{\mu(t)-\mu(\delta+t)}{\mu(4)}
$$ for every $\delta \in (0,1/2)$ and $t \in [t_*, 2]$, by using the continuity of $\mu$, we choose $\delta_2 \in (0, \delta_1]$ such that
\bee \label{e:mu04}
  \mu(t)  \le (1+\eps)\, \mu\big(t(1-\delta_2)^{-1}\big) ~  \text{ and } ~\mu(t)  \le (1+\eps)\, \mu(t+\delta_2)     \quad \text{for every } t \in [t_*, 2].
 \eee
Combining \eqref{e:mu02}--\eqref{e:mu04}, we have that for every $\delta \le \delta_2$,
\bee \label{e:mu05}
 \mu(t)  \le (1+\eps)^{4+\alpha/2} \times
 \begin{cases}
  \mu\big(t(1-\delta)^{-1}\big) \quad &\text{ when } t <2\\
   \mu(t+\delta)\quad &\text{ when } t\ge 1/2.
 \end{cases}
 \eee

Let $r>2$. Using \eqref{j11}, we
put
$$\ds j(r+\dt) =\left( \int_0^1 + \int_1^{\ift} \right) (4\pi t)^{-d/2} \exp \big(-\frac{(r+\dt)^2}{4t}\big) \, \mu (t)\, dt =: I + II. $$
Since  $(1-\dt)(r+\dt)^2  \le r^2 + \dt (r+\dt)\big( 2- (r+\dt) \big) \le r^2$, by \eqref{e:mu05} and a change of variables,
\begin{eqnarray*}
I &\geq& \int_0^1 (4\pi t)^{-d/2} \exp \big(-\frac{(1-\dt)^{-1} r^2}{4t}\big) \, \mu (t) \,dt \\
&=& (1-\dt)^{-1+d/2} \int_0^{1-\dt} (4\pi t)^{-d/2} \exp (-\frac{r^2}{4t}) \,\mu \big(t(1-\dt)^{-1}\big) \,dt \\
&\geq& (1-\dt)^{-1+d/2} (1+\ep)^{-4-\af/2} \int_0^{1-\dt} (4\pi t)^{-d/2} \exp (-\frac{r^2}{4t}) \,\mu(t) \,dt \quad \text{for every}\,\, \delta \le \delta_2 .
\end{eqnarray*}

On the other hand,
from $0 \leq (r+\dt-t)^2  = (r+\dt)^2 -2tr + t(t-\dt)- \dt t$,
we see that
$ t(t-\dt) \geq 2tr + \dt t - (r+\dt)^2$. Thus  we get
\begin{eqnarray*}
\frac{(r+\dt)^2}{4t} - \frac{r^2}{4(t-\dt)} = \frac{(r+\dt)^2 (t-\dt) - r^2 t}{4t(t-\dt)}
= \frac{\dt (2tr+\dt t -(r+\dt)^2)}{4t(t-\dt)} \leq \frac{\dt}{4} .
\end{eqnarray*}
Therefore by using this, a change of variables, \eqref{e:mu05} and the inequality $  t+\dt  \le t \,(1-\dt)^{-1}  $
for $1-\dt \le t<\ift $, we get
\begin{eqnarray*}
II &\geq& e^{-\dt /4} \int_1^{\ift} (4\pi t)^{-d/2} \exp\big( -\frac{r^2}{4(t-\dt)}\big) \,\mu (t) \,dt \\
&=& e^{-\dt /4} \int_{1-\dt}^{\ift} (4\pi (t+\dt))^{-d/2} \exp(-\frac{r^2}{4t}) \,\mu(t+\dt) \,dt\\
&\geq& e^{-\dt /4} (1+\ep)^{-4-\af/2} (1-\dt)^{d/2} \int_{1-\dt}^{\ift} (4\pi t)^{-d/2} \exp(-\frac{r^2}{4t}) \,\mu(t)\,dt \quad \text{for every} \,\, \delta \le \delta_2 \,.
\end{eqnarray*}

Consequently for every $\delta \le \delta_2$ and $ r>2$,
$$ j(r+\dt) \geq \big( (1-\dt)^{-1+d/2} \,\wedge\, e^{-\dt /4} (1-\dt)^{d/2} \big) (1+\ep)^{-4-\af/2} j(r)  $$
and so
$$
\limsup_{\delta \downarrow 0} \left(\sup_{r>2}\frac{j(r)}{j(r+\delta)}\right) \le (1+\ep)^{4+\af/2}.
$$Since $\ep>0$ is arbitrary and $\frac{j(r)}{j(r+\delta)} \ge 1$, the proof is completed.
\qed

\begin{lemma}\label{l:j2}
\begin{equation*}
\lim_{\delta \downarrow 0} \sup_{r \in (0, 4]}\frac{j(r)}{j\big(r(1+\delta)\big)} = 1\,.
\end{equation*}
\end{lemma}
\pf
Fix $\eps>0$ and let
$\sA:=\alpha\Gamma((d+\alpha)/2){2^{-1+\alpha}\pi^{-d/2}(\Gamma(1-
\alpha/2))^{-1}}$.
By Potter's Theorem \cite[Theorem 1.5.6(i)]{BGT},
there exists $r_1=r_1(\eps)>0$
such that
$$
\frac{\ell(t^{-2})}{\ell(s^{-2} )} \,\ge\,  (1+\eps)^{-1} \min \left\{\frac{t}{s} \,,\,
\frac{s}{t} \right\}\quad \text{ for }s, t\le 2\,r_1.
$$
Moreover by Theorem \ref{t:Gorigin}, there exists $r_2=r_2(\eps)>0$
such that
$$
1+\eps \,\ge\, \frac{ \sA \ell (s^{-2} )     }{s^{d+\alpha} j(s)} \,\ge\,  (1+\eps)^{-1}\quad  \text{for }s \le 2 \, r_2.
$$
Thus for $r \le r_3:=r_1 \wedge r_2$ and $\delta \in (0,1)$
\begin{eqnarray*}
\frac{j\big(r(1+\delta)\big)}{j(r)}&=&
\left(\frac{j\big(r(1+\delta)\big) \, r^{d+\alpha} \,(1+\delta)^{d+\alpha}}{\sA\, \ell \big(r^{-2} (1+\delta)^{-2}\big)}
\right) \left(\frac{ \sA \, \ell (r^{-2} )     }{r^{d+\alpha} j(r)}\right)
\frac{ \ell \big(r^{-2} (1+\delta)^{-2}\big)}{  \ell (r^{-2} ) } \, (1+\delta)^{-d-\alpha}\\
& \ge& (1+\eps)^{-3}(1+\delta)^{-d-\alpha-1} \, .
\end{eqnarray*}
On the other hand  for every $\delta \in (0,1)$ and $r \in [r_3, 4]$,
$$\frac{j(r)-j\big((1+\delta)r\big)}{j\big((1+\delta)r\big)} \,\le\, \frac{j(r)-j\big((1+\delta)r\big)}{j(8)} \,\le\, j(8)^{-1} \delta r  | j\,'\,(r_3)| \le 4j(8)^{-1} \dt | j\,'\,(r_3)| $$
and so
$
\big( 1+4j(8)^{-1} \dt  | j \, ' \, (r_3)| \big) \, j(r(1+\delta))
 \,\ge\,  j(r)
$. Therefore
$$
\limsup_{\delta \downarrow 0} \left( \sup_{r \in (0, 4]}\frac{j
(r
)}{j
\big(r
(1+\delta)\big)}  \right) \le (1+ \eps)^3.
$$
Since $\eps>0$ is arbitrary and $\frac{j
(r
)}{j
(r
(1+\delta))} \ge 1$, we complete the proof.
\qed

In this section, for
 the notational convention we define$$
\Lambda_{a, b}({u}):= \int_{A(0, a, b)} j(|y|) {u}(y)dy  \quad \text{and} \quad \Lambda_{a}({u}):= \int_{B_a ^{\, c}} j(|y|) {u}(y)dy
$$ for every nonnegative function $u$ on $\R^d$ and constants $a$ and $b$ with $b>a>0$. By Lemmas \ref{l:j1} and \ref{l:j2}, there exists an increasing continuous function $\delta (\eps) : (0,1/2] \to (0,1/2]$ such that $\lim_{\eps \downarrow 0} \delta (\eps)=0$ and
\bee\label{e:jeps}
 \left(\sup_{r
  > 2}\frac{j(r)}{j(r+\delta(\eps))}\right) \vee  \left(\sup_{r \in (0, 4]}\frac{j(r)}{j(r(1+\delta(\eps) ))} \right) \, \le\, 1+ \eps.
\eee

\begin{lemma}  \label{l:p1}
For every $0< \eps \le1/2$, $0<p \le 1/2$ , $r \le 2$ and any nonnegative function $u$ in $\RR^d$, we have for every $x \in B_{\dt  pr /3} $
$$
  (1 + \varepsilon)^{-1}\Lambda_{ pr}(u) \E_x[\tau_{B_{\dt  pr /3}}]
  \,\le\, \int_{B_{pr}^c }P_{B_{\dt  pr /3}}(x,y)u(y)dy
 \, \le \,  (1 + \varepsilon) \Lambda_{ pr}(u) \E_x[\tau_{B_{\dt  pr /3}}]
$$
where $\delta=\delta(\eps) \in (0,1/2]$ is in \eqref{e:jeps}.
\end{lemma}
\pf
If  $z \in B_{\dt  pr /3}$ and $y \in A(0, pr, 1)$, then we have
$$|y-z| \,\le\, |y|+|z| \,\le\, |y| + \dt  pr /3 \,\le\, (1+\dt/3)|y| \,\le\, (1+\delta
)|y|$$
 $$\text{and} \quad |y-z| \,\ge\, |y|-|z| \,\ge\, |y| - \dt  pr /3 \,\ge\, (1-\dt/3)|y| \,\ge\, (1+\delta
)^{-1}|y| .$$
Thus by \eqref{e:jeps} and the fact that $r \mapsto j(r)$ is decreasing,
\begin{eqnarray*}
1+\eps \,\ge\, \frac{j( (1+\delta
)^{-1}|y|)}{j(|y|)} \,\ge\, \frac{j(|y-z|)}{j(|y|)} \,\ge\, \frac{j( (1+\delta
)|y|)}{j(|y|)}
\,\ge\, (1+\eps)^{-1}  \quad \text{for } y \in A(0, pr,  1) \,.
\end{eqnarray*}
On the other hand, since the assumptions $r \le 2$ and $p \le 1/2$ imply
$
\dt  pr /3 \le \delta
,
$
we have $$|y-z| \le |y|+|z| \le |y| + \dt  pr /3 \le   |y| + \delta
    $$
 $$\text{and} \quad |y-z| \ge |y|-|z| \ge |y| - \dt  pr /3 \ge |y| - \delta .
  $$ Thus
by \eqref{e:jeps} and the fact that $j$ is decreasing,
\begin{eqnarray*}
1+\eps \,\ge\, \frac{j ( |y| - \delta
 )}{j(|y|)}\,\ge\, \frac{j(|y-z|)}{j(|y|)} \,\ge\, \frac{j(  |y| + \delta
 )}{j(|y|)}
\,\ge\, (1+\eps)^{-1} \quad\text{for } |y| \ge 1 \,.
\end{eqnarray*}

So we have  for $x \in B_{\dt  pr /3}$,
\begin{eqnarray*}
&&   \int_{B_{pr}^c }P_{B_{{\dt  pr /3}}}(x,y)u(y)dy  =
  \int_{B_{pr}^c} \int_{B_{{\dt  pr /3}}} G_{B_{{\dt  pr /3}}}(x, z) j(|z- y|) dz \, u(y)dy\\
 &&  \leq
  (1 + \varepsilon)
   \int_{B_{{\dt  pr /3}}} G_{B_{{\dt  pr /3}}}(x, z) dz  \int_{B_{pr}^c}j(|y|)  u(y)dy\,= \,
    (1+\eps) \E_x[\tau_{B_{{\dt  pr /3}}}]
   \Lambda_{ pr}(u)
\end{eqnarray*}
and
\begin{eqnarray*}
  \int_{B_{pr}^c }P_{B_{{\dt  pr /3}}}(x,y)u(y)dy
  &\ge&
 (1 + \varepsilon)^{-1} \int_{B_{\dt  pr /3}} G_{B_{\dt  pr /3}}(x, z) dz  \int_{B_{pr}^c}j(|y|)u(y)dy \, \\
 &=&
  (1 + \varepsilon)^{-1}  \E_x[\tau_{B_{\dt  pr /3}}] \Lambda_{ pr}(u).
\end{eqnarray*}
\qed

The next two results were proved in \cite{KSV4} in a more general setting.
\begin{lemma}{\rm(\cite[Lemma 5.2]{KSV4})}\label{l:p2}
For every $p \in (0,1)$, there exists $c=c(\alpha, d,\ell, p)>0$ such that for every
$r \in (0, 1)$ and  $ (x,y) \in  B_{pr} \times B_r^c$,
\begin{eqnarray*}
P_{{B_r}}(x, y)
\le
\,\frac{c}{\phi(r^{-2})}\left(
\int_{A(0, (1+p)r/2 ,r)}  j(|z|)  P_{{B_r}}(z,y)dz +j(|y|)
\right).\end{eqnarray*}
\end{lemma}
\begin{lemma}{\rm(\cite[Lemma 5.4]{KSV4})}\label{l:p3}
 There exists $c=c(  \alpha,d, \ell)>1$ such that
 for every $r \in (0, 1)$ and $ (x,y) \in B_ {r/2} \times B_r^c$,
$$P_{B_r}(x,y) \ge  \frac{c}{\phi(r^{-2})} \left(\int_{A(0, r/2, r) }  j(|z|)  P_{B_r}(z,y)dz  +j(|y|) \right) .\\
$$
\end{lemma}
Note that since $\ell$ is slowly varying at $\infty$ and $\ell$ is strictly positive and continuous on $(0, \infty)$, there exists a
constant $c=c( \alpha, \ell) >1 $ such that for every
$r \in (0, 1)$,
\begin{equation}\label{lll}
c^{-1} \,\le\,\frac{\ell\big((2 r/3)^{-2}\big)}{\ell(r^{-2})} \,\le\,  \left(\frac{\ell\big((2 r/3)^{-2}\big)}{\ell(r^{-2})} \vee
\frac{\ell \big((r/2)^{-2}\big)}{\ell(r^{-2})}
\right) \,\le\, c.
\end{equation}

Recall that $C_0$ is the constant in Theorem \ref{HP}.
\begin{lemma}\label{l:u}
There exists $C_*=C_*(\alpha, d, \ell) \,\ge\, C_0$ such that for every $r \in (0, 1)$, any nonnegative function $u$ in $\RR^d$
which is regular harmonic in ${B_r}$ with respect to $X$ and for any  $ x \in  B_{r/2}$,
\begin{align}
C_*^{-1} \E_x[\tau_{B_r}]
\Lambda_{r/2} (u) \,\le\, u(x)\,&\,\le\, C_*\, \E_x[\tau_{B_{2r/3}}]\Lambda_{3r/4} (u)\label{*1}\\
 &\,\le\, C_*\, \E_x[\tau_{B_r}]\Lambda_{r/2} (u).\label{*2}
\end{align}
\end{lemma}
\pf
Since  $u$ is regular harmonic in ${B_r}$ with respect to $X$ and $\P_z (X_{\tau_{B_r}} \in \partial B_r) = 0$ for $ z \in B_r$, we have
$u(z)= \int_{B_r ^c}  P_{B_r}(z,y) u(y) dy $ for every $z \in B_r$ (see \eqref{e:harmonic}).
Thus by using Lemma \ref{l:p2} in the first, and \eqref{lll} in the second inequality, we get
\begin{eqnarray*}
u(x) &\le & \frac{c_1}{\phi(r^{-2})} \left(\int_{B_r ^c}\int_{A(0, 3r/4, r)}  j(|z|)  P_{B_r}(z,y)dz u(y)dy + \int_{B_r^c} j(|y|) u(y)dy\right)\\
&= & \frac{c_1}{\phi(r^{-2})} \left(\int_{A(0, 3r/4, r)}  j(|z|) \left( \int_{B_r^c}  P_{B_r}(z,y) u(y)dy \right) dz + \int_{B_r^c} j(|y|) u(y)dy\right) \\
&= & \frac{c_1}{\phi(r^{-2})} \left(\int_{A(0, 3r/4, r)}  j(|z|) u(z) dz+\int_{B_r^c} j(|y|) u(y)dy\right)\\
& \le & \frac{c_2}{\phi\big((2r/3)^{-2}\big)} \int_{B_{3r/4}^{\, c}} j(|y|) u(y)dy.
 \end{eqnarray*}
Similarly using Lemma \ref{l:p3}, we also get
$ u(x) \ge \frac{c_3}{\phi(r^{-2})} \int_{B_{r/2}^{\, c}} j(|y|) u(y)dy.$
Now applying \cite[Lemmas 13.4.2 and 13.4.3]{KSV2}, we have proved \eqref{*1}. \eqref{*2} follows immediately from \eqref{*1}. \qed

For the remainder of the section, we fix $C_*$ in Lemma \ref{l:u}.
\begin{lemma} \label{l:osc}
Suppose that $r \in (0,1)$.
For nonnegative functions $u_1, u_2$ in $\RR^d$
which are harmonic in $B_r$ with respect to $X$, we have
for every $0 < p < q/4<  1/8$,
\begin{equation*}
  \expr{\sup_{B_ {pr}} \frac{g_1}{g_2} -
    \inf_{B_{pr}} \frac{g_1}{g_2}} \leq \frac{
  C_*^2 - 1}{C_*^2 + 1} \left(\sup_{B_{qr}} \frac{u_1}{ u_2} - \inf_{B_{qr}} \frac{u_1}{ u_2}\right),
\end{equation*}
where
$g_i(x):=\E_x[u_i(X_{\tau_{B_{2pr}}}): X_{\tau_{ B_{2pr}}}\in A(0, 2  pr ,  q r)  ].$\end{lemma}
\pf
For $a >0$, we define $m_a = \inf_{B_{a}} (u_1 / u_2)$ and $M_a
= \sup_{B_{a}} (u_1 / u_2)$.
Let   $$
  f(x):= \E_x[(u_1-m_{qr}u_2)(X_{\tau_{B_{ 2p r}}}): X_{\tau_{B_{2 p r}}}\in A(0, 2  pr ,  q r)  ]  =
  g_1 (x)- m_{{qr}} g_2 (x)
  $$
  and
   $$
  h(x):= \E_x[(M_{qr}  u_2- u_1)(X_{\tau_{B_{ 2p r}}}): X_{\tau_{B_{2 p r}}}\in A(0, 2  pr ,  q r)  ] =
M_{qr}  g_2 (x)-  g_1 (x) , $$
then $f$ and $h$ are regular harmonic in $B_{2pr}$ and nonnegative in $\R^d$.
Thus by 
applying \eqref{*2} to $f$ and $h$, we get
 \begin{eqnarray*}
  \sup_{B_{{pr}}} \frac{g_1}{g_2} - m_{qr} \,=\,
  \sup_{B_{{pr}}} \frac{f}{g_2} \,\leq\,
  C_*^2 \inf_{B_{{pr}}} \frac{f}{g_2} \,=\,
  C_*^2 \expr{\inf_{B_{{pr}}} \frac{g_1}{g_2} - m_{{qr}}}
\end{eqnarray*}
and
\begin{eqnarray*}
  M_{qr} - \inf_{B_{{pr}}} \frac{g_1}{g_2} \,=\,
  \sup_{B_{{pr}}} \frac{h}{g_2} \,\leq\,
  C_*^2 \inf_{B_{{pr}}} \frac{h}{g_2} \,=\,
  C_*^2 \expr{M_{qr} - \sup_{B_{{pr}}} \frac{g_1}{g_2}} \,.
\end{eqnarray*}
By adding these inequalities, we proved the lemma.
\qed

Now we are ready  prove the main result of this section. We prove the main result for the quotient of two harmonic functions  in the next theorem.
We closely follow the proof of \cite[Lemma 8]{BKK}.

\begin{thm}\label{th:osc_1}
For every $\eta > 0$, there exists $a=a(\eta, \af, d, \ell) \in (0, 1)$  such that for every $x_0 \in \R^d$ and $r \in (0, 1]$,
\begin{equation*}
\sup_{ B(x_0,ar)} \frac{u_1}{u_2} \leq (1 + \eta) \inf_{ B(x_0,ar)} \frac{u_1}{u_2}
\end{equation*}
for every nonnegative functions $u_1$ and $u_2$ in $\RR^d$
which are harmonic in $B(x_0, r)$ with respect to $X$.
\end{thm}
\pf
We assume $x_0=0$. We fix $r \in (0, 1]$ and  nonnegative functions $u_1, u_2$ in $\RR^d$
which are 
harmonic in $B_r$ with respect to $X$. Fix ${\eta} >0$ and let $$\ds \varphi(t) := 1 + \frac{\et}{2(C_* ^2 +1)} + \frac{C_*^2}{C_*^2 + 1} (t - 1) \quad \text{for }t\geq 1\, \quad\text{and}\quad \varphi^{1}\,:=\,\varphi, \,\, \varphi^{l+1}\,:=\,\varphi(\varphi^{l}) \,\,\text{ for  } l=1,2,\cdots.$$
Then $$\ds \varphi^l(C_*^2) \,=\, 1+ \frac{\et}{2(C_* ^2 +1)} \sum_{i=0}^{l-1} (\frac{C_*^2}{C_*^2 + 1})^i + (\frac{C_*^2}{C_*^2 + 1})^l (C_*^2-1) \,\le\,
1+ \frac{{\eta}}{2} + (\frac{C_*^2}{C_*^2 + 1})^l (C_*^2-1).$$
Choose
 $l= l(C_*, {\eta})$ large such that \begin{equation} \label{e:Ceta4}(\frac{C_*^2}{C_*^2 + 1})^l \,(C_*^2-1) <\frac{{\eta}}{2} \quad \text{ so that } \quad   \varphi^l(C_*^2)<1+\eta\,.\end{equation}
Also we choose $\varepsilon = \varepsilon(\eta)$ small enough so that
\begin{equation} \label{e:Ceta}
1+\frac{\et}{C_* ^2 +1} \ge
 \big(C_*^3 \, \varepsilon + (1 + \varepsilon)\big)^2 (1 + \varepsilon)^{2 } ,
\end{equation}
\begin{equation} \label{e:Ceta2}(1 + C_*^2 \varepsilon)^2
 \le 1 + \frac{\et}{2(C_* ^2 +1)}\quad \text{and} \quad 1+ C_*^2 \, \varepsilon  \le \frac{C_*^2}{C_*^2 - 1}  . \end{equation}
Let $k= k(\varepsilon)\ge 3$ be the smallest integer such that $k  > 1+
1/\varepsilon^{2}$.
We recall that $\delta=\delta(\eps)>0$ is the constant from \eqref{e:jeps} and fix it. Let $\ds p_i :=
(\delta /6)^i/2$ 
for $i=0,\cdots,lk-1$.
For simplicity, 
we put $m_a := \inf_{B_{a}} u_1 / u_2$ and $M_a
:= \sup_{B_{a}} u_1 / u_2$.

\medskip
{\it Case 1.}
Suppose that  the following holds for both $ i = 1\text{ and
} 2$; for every $ 0 \le \, m<\, lk$,
\begin{equation*}
\int_{A(0, rp_{m+1}, rp_m)} j(|y|) u_i (y)\, dy \,=\, \Lambda_{r p_{m +
    1}, rp_{m}}(u_i) \, >\, \varepsilon \Lambda_{r p_{m}}(u_i)  \,=\, \varepsilon \int_{B_{rp_m} ^c} j(|y|) u_i (y) \,dy\,.
\end{equation*}

By the definition
of $k$, for $0\leq j \le l-1$
\begin{eqnarray}
\Lambda_{2 r p_{(j + 1) k},r p_{j k}}(u_i)&\ge& \Lambda_{r p_{(j +1)  k - 1},r p_{j k}}(u_i) \,
  =\,\sum_{m=0}^{k-2}  \Lambda_{r p_{jk+m+1}, rp_{j k+m}}(u_i) \nonumber \\
&\ge&  \varepsilon  \sum_{m=0}^{k-2}    \Lambda_{r p_{j k+m}}(u_i) \,\ge\,   (k - 1)  \varepsilon \, \Lambda_{r p_{j k}}(u_i) \,\geq\,
  \varepsilon^{-1} \Lambda_{r p_{j k}}(u_i) \,. \label{e:Leps}  
\end{eqnarray}

For $i= 1,\,2$ and $j = 1, \cdots, l-1$, we let
\begin{equation*} 
f_i^j(x) \,:= \, \E_x[u_i (X_{\tau_{B_{2r p_{(j + 1) k}}}}) : X_{\tau_{B_{2r p_{(j + 1) k}}}} \in B_{rp_{jk}}^{\,c}] \,=\, \int_{ B_{rp_{jk}}^{\,c}} P_{B_{2r p_{(j + 1) k}}} (x,y) \, u_i (y) \, dy 
\end{equation*}
and
\begin{eqnarray*} 
 g_i^j (x)&:=& \E_x[u_i (X_{\tau_{B_{2r p_{(j + 1) k}}}}) : X_{\tau_{B_{2r p_{(j + 1) k}}}} \in A(0, 2rp_{(j+1)k} , rp_{jk})]  \\
&=& \int_{  A(0, 2rp_{(j+1)k} , rp_{jk})} P_{B_{2r p_{(j + 1) k}}} (x,y) \, u_i (y) \, dy\, ,
\end{eqnarray*}
which  are regular harmonic in $B_{2rp_{(j+1)k}}$ and $u_i = f_i ^j +g_i ^j$.

By \eqref{*1} applied to $B_{rp_{(j + 1) k}} $ in the first, and the facts that $f_i^j (x) =0$ on $A(0, 2 r p_{(j + 1) k}, r p_{j k})$ and $f_i^j (x)=u_i(x)$ on
 $B_{rp_{jk}}^{\,c}$ in the second inequality,
 we have for $x \in B_{rp_{(j+1)k}}$,
 $$f_i^j(x)
\leq C_{*} \, \E_x[\tau_{B_{\frac{4}{3}r p_{(j + 1) k}}}]\, \Lambda_{\frac{3}{2} rp_{(j+1)k}} (f_i^j)
\, \le \,  C_{*} \,\E_x[\tau_{B_{2r p_{(j + 1) k}}}]\, \Lambda_{r p_{j k}}(u_i)   \quad\text{for }\,  j= 1, \cdots, l-1.
$$
Hence by \eqref{e:Leps}, the fact that $g_i^j(x) =u_i(x)$ on $A(0, 2  p_{(j + 1) k}r,  p_{j k} r)$ and
\eqref{*2} applied to
$B_{rp_{(j + 1) k}}$,
\begin{eqnarray*}
f_i^j(x) &\leq&
  C_{*}\, \varepsilon \, \E_x[\tau_{B_{2r p_{(j + 1) k}}}]\, \Lambda_{ 2 rp_{(j + 1) k}, rp_{j k}}(u_i)
 \,=\,
  C_{*}\, \varepsilon \, \E_x[\tau_{B_{2r p_{(j + 1) k}}}]\, \Lambda_{ 2r p_{(j + 1) k}}(g_i^j)
\\
&\le&
  C_{*}\, \varepsilon \, \E_x[\tau_{B_{2r p_{(j + 1) k}}}]\, \Lambda_{ r p_{(j + 1) k}}(g_i^j)
\,\leq\, C_{*}^2\, \varepsilon\,  g_i^j(x) \quad \text{for} \,\, x \in B_{rp_{(j + 1) k}} \text{ and } j= 1, \cdots , l-1.
\end{eqnarray*}
Since $u_i (x)=
f_i^j(x)+g_i^j(x)$  and $\ds \frac{g_1^j}{f_2^j + g_2^j} \le  \frac{u_1}{u_2}  \le \frac{f_1^j + g_1^j}{ g_2^j}  $, we have  \begin{eqnarray*}
(1 + C_{*}^2 \, \varepsilon)^{-1} \inf_{B_{rp_{(j + 1) k}}}
  \frac{g_1^j}{g_2^j}  \le m_{rp_{(j + 1) k}} \le M_{rp_{(j + 1) k}}
   \le (1 + C_{*}^2 \, \varepsilon) \sup_{B_{rp_{(j + 1) k}}}\frac{g_1^j}{g_2^j} \, , \quad  j= 1, \cdots, l-1 \, .
\end{eqnarray*}
Thus by Lemma~\ref{l:osc},
\begin{eqnarray*}
&&(C_*^2 + 1) \left((1 + C_{*}^2 \varepsilon)^{-1}M_{rp_{(j + 1) k}}-(1 + C_{*}^2 \varepsilon)
m_{rp_{(j + 1) k}}\right)\\
&&\leq  (C_*^2 + 1)\left(\sup_{B_{rp_{(j + 1) k}}} \frac{g_1^j}{g_2^j}  -
    \inf_{B_{rp_{(j + 1) k}}} \frac{g_1^j}{g_2^j} \right)\leq
(C_*^2 - 1) (M_{rp_{j k}} - m_{rp_{j k}})\, , \quad  j= 1, \cdots, l-1\,.
\end{eqnarray*}
Multiplying by $(1 + C_{*}^2\, \varepsilon)/(m_{rp_{(j + 1) k}}(C_*^2 + 1))$ and using the obvious fact $m_{rp_{(j + 1) k}} \geq m_{rp_{j k}}$, we obtain
\begin{equation*}
  \frac{M_{rp_{(j + 1) k}}}{m_{rp_{(j + 1) k}}}
 \le  (1 + C_{*}^2 \varepsilon)^2 +
    (1 + C_{*}^2 \varepsilon)
\frac{C_*^2 - 1}{C_*^2 + 1}
 \left(  \frac{M_{rp_{j  k}}}{m_{rp_{j  k}}} - 1  \right) \,.
\end{equation*}
By the definition of $\varphi$ and \eqref{e:Ceta2},  $\ds
  \frac{M_{rp_{(j+1)k}}}{m_{rp_{(j+1)k}}}
  \leq
  \varphi  \left(
  \frac{M_{rp_{jk}}}{m_{rp_{jk}}}
   \right) .$ We already know that  $\frac{M_{r/2}}{m_{r/2}} \leq C_*^{\,2}$ by 
\eqref{*2}.
And also by the monotonicity of $\varphi$ and \eqref{e:Ceta4}, we get $$\ds
  \frac{M_{rp_{lk}}}{m_{rp_{lk}}}
  \,\leq\,
  \varphi  \left(
  \frac{M_{rp_{(l - 1)k}}}{m_{rp_{(l - 1)k}}}
   \right)\,\leq\,
  \cdots \,\leq\,
  \varphi^l   \left(\frac{M_{r/2}}{m_{r/2}} \right) \,\leq\, \varphi^l (C_*^{\,2}) \,<\,
  1 + {\eta}  \, .
$$

\noindent
\medskip

{\it Case 2.}
Suppose that there exists $m < lk$ such that for either $i=1$ or $2$,
\begin{equation*} 
   \int_{A(0, r p_{m +1} ,r p_{m} )}j(|y|)u_i(y)\,dy \,=\, \Lambda_{r p_{m +1}, r p_{m}}(u_i) \,\leq\,
  \varepsilon \, \Lambda_{r p_{m}}(u_i)\,=\,\varepsilon \int_{B_{rp_{m}} ^c}j(|y|)u_i(y)\,dy \,.
\end{equation*}

Note that by 
\eqref{*2}, 
$$\ds C_{*}^{-1}\frac{u_{3-i}(y)}{\Lambda_{r p_{m}}(u_{3-i})} \,\leq\,
  \E_y[\tau
 _{B_{2r p_{m}}}] \,\leq\,
  C_{*} \frac{u_{i}(y)}{\Lambda_{r p_{m}}(u_{i})} \quad \text{for }  y \in A(0,rp_{m + 1}, r p_{m} ) .$$
Hence by integrating on $ A(0,rp_{m + 1}, r p_{m} )$, we get
$$ \frac{\Lambda_{r p_{m+1},rp_{m}}(u_{3-i})}
{\Lambda_{rp_{m}}(u_{3-i})} \,\leq\,
C_{*}^2 \frac{\Lambda_{r p_{m+1},rp_{m}}(u_{i})}
{\Lambda_{rp_{m}}(u_{i})}\, \le C_{*}^2 \varepsilon \,.$$ 
%Thus there exists $m <lk$ such that
Thus
\bee \label{e:assume1}
 \Lambda_{r p_{m + 1}, rp_{m}}(u_i) \le C_{*}^{2} \varepsilon \,
\Lambda_{r p_{m}}(u_i)\quad \text{for both } i = 1 \text{ and } 2\,.
\eee

%For such $m$, let
Let 
\begin{equation*}
f^{m}_i(x)= f_i(x)\,:=\,\E_x[u_i ( X_{\tau_{B_{2r  p_{m + 1 }}}} ): X_{\tau_{B_{2r  p_{m + 1 }}}} \in {B_{rp_{m}} ^c}]  \,=\, \int_{B_{rp_m} ^c} P_{B_{2rp_{m+1}}} (x,y) \, u_i (y) \,dy   
\end{equation*}
and 
\begin{eqnarray*} g^{m}_i(x) = g_i(x)&:=&\E_x[u_i (X_{\tau_{B_{2r  p_{m + 1 }}}}): X_{\tau_{B_{2r  p_{m + 1 }}}} \in A(0, 2  rp_{m + 1 } ,  r p_{m } )  ] \\
&=& \int_{A(0, 2rp_{m+1}, rp_m)} P_{B_{2rp_{m+1}}} (x,y) \, u_i (y) \,dy \, , \end{eqnarray*}
so that $u_i = f_i +g_i $.
Since $g_i$ is regular harmonic in $B_{2rp_{m+1}}$, by \eqref{*1} %Lemma \ref{l:u} 
we obtain for $x \in B_{r p_{m+1 } }$,
\begin{equation*}
g_i(x)  \,\leq\,   C_* \, \E_x[\tau_{B_{\frac{4}{3} r p_{m+1 } }}]
 \Lambda_{\frac{3}{2}rp_{m+1}} (g_i) \,\le\, %\,=\,
  C_{*} \, \E_x[\tau_{B_{2 r p_{m+1 } }}]
 \Lambda_{rp_{m+1}} (g_i) \,.
\end{equation*}
Also since $g_i=0$ on $\overline{B_{rp_{m}} }^c$ and $g_i=u_i$ on $A(0, 2r  p_{m + 1 },  r p_{m} )$, we get
\begin{eqnarray*}
g_i(x) & \leq &  C_* \, \E_x[\tau_{B_{2r p_{m+1 } }}]
 \Lambda_{rp_{m+1}, rp_m } (g_i) \,\le\,
   C_* \, \E_x[\tau_{B_{2r p_{m+1 } }}]
 \Lambda_{rp_{m+1}, rp_m } (u_i)  \\
&\leq& \varepsilon \,C_* ^3 \,\E_x[\tau_{B_{2r p_{m+1 } }}] \Lambda_{rp_{m+1}} (u_i) \quad \text{for}\quad x \in B_{r p_{m+1 } } .
  \end{eqnarray*}
The last inequality comes from \eqref{e:assume1}.

Then by \eqref{e:assume1}, applying Lemma \ref{l:p1} to $f_i (x)$ and the fact that $ \ds \frac{f_1}{f_2 + g_2} \le \frac{u_1}{u_2} \le \frac{f_1 + g_1}{f_2} $, %on $B_{2rp_{m+1}}$,
we have
$$   \frac{(1 + \varepsilon)^{-1} \Lambda_{r p_{m}}(u_1)}
    {\big(   (1 + \varepsilon) + \varepsilon \, C_*^3\big) \Lambda_{r p_{m}}(u_2)} \,\leq\,
 \frac{u_1 (x)}{u_2 (x)} \,\leq\,
  \frac{\big(  (1 + \varepsilon) + \varepsilon \, C_*^3 \big) \Lambda_{r p_{m}}(u_1)}
    {(1 + \varepsilon)^{-1} \Lambda_{r p_{m}}(u_2)}
  \,  \quad \text{ for } x \in B_{rp_{m + 1 }}.
$$So by \eqref{e:Ceta}, $\ds
  \frac{ M_{rp_{lk }}}
  {m_{rp_{lk }}} \,\leq\,
  \frac{ M_{rp_{m + 1 }}}
  {m_{rp_{m + 1 }}} \,\le\, \big( \varepsilon \, C_*^3 + (1 + \varepsilon)\big)^2 (1 + \varepsilon)^{2 } \,\leq\,
  1 + \frac{\et}{C_* ^2 +1} \,<\, 1+\et
$ .

In these two cases, we prove the theorem with $a= p_{lk}$.
\qed

\noindent
{\bf Proof of Theorem \ref{th:osc}.}
Take $u_1 = u$ and $u_2 \equiv 1$ in Theorem \ref{th:osc_1}.
\qed

As a corollary of Theorem \ref{th:osc}, we get
\begin{cor}\label{c:conti}
There exists an increasing continuous function $\theta:(0, 1) \to (0, \infty)$ with $\lim_{t \to 0} \theta(t)=0$ such that
 for every $x_0 \in \R^d$,  $R \in (0, 1]$ and $r < R/2$,
$$
\sup_{x,y \in B(x_0,R/2), \, |x-y|<r} |u(x)-u(y)|   \leq
\theta (|x-y|/r) \sup_{w \in B(x_0, R)} |u(w)|
$$for nonnegative function $u$ in $\RR^d$
which is harmonic in $B(x_0, R)$ with respect to $X$.
\end{cor}
\pf
Without loss of generality, we assume $x_0 = 0$. For fixed $R \in (0,1]$ and $r$ with $r<R/2$, let $x,y \in B_{R/2}$ be such that $|x-y|<r$ and $x, y \in B(z, |x-y|) \subset B_{R}$ for some $z \in B_{R/2}$.
For a nonnegative integer $k$, by Theorem \ref{th:osc} we can choose $a_{k+1} < a_k$ recurrently such that
\begin{equation} \label{e:osc_cor}
 \sup_{B(z, r a_k)} u \,\le\, (1+2^{-k-1}) \inf_{B(z, r a_k)} u \quad \text{for  } z \in B_{R/2}.
\end{equation}

Define $a(\eta)$ using the linear interpolation as
\begin{displaymath}
a(\eta) = \left\{ \begin{array}{ll}
a_k & \text{if} \quad \eta = 2^{-k} \\
\ds \frac{a_k - a_{k+1}}{2^{-k} - 2^{-k-1}} \, \eta +2 a_{k+1} - a_k & \text{if} \quad 2^{-k-1} < \eta < 2^{-k} .
\end{array} \right.
\end{displaymath}
Then $a(\eta)$ is continuous and strictly increasing, so there exists an inverse function $\theta := a^{-1} : (0,1) \rightarrow (0, \ift)$, which is increasing and continuous.

Now we choose a nonnegative integer $k$ such that $\ds a_{k+1} \le \frac{|x-y|}{r} < a_k$, so that $\ds 2^{-k-1} \le \theta \big( \frac{|x-y|}{r} \big)$. Using this and \eqref{e:osc_cor}, we get
\begin{eqnarray*}\sup_{B(z,|x-y|)} u &\le& \sup_{B(z,ra_k)} u \,\le\, (1+2^{-k-1}) \inf_{B(z,ra_k)} u \,\le\, \big(1 \,+\, \theta(\frac{|x-y|}{r}) \big) \inf_{B(z,ra_k)} u \\
&\le& \big(1\,+\, \theta(\frac{|x-y|}{r}) \big) \inf_{B(z,|x-y|)} u \,.
\end{eqnarray*}
Therefore
$$
|u(x)- u(y)| \,\le\, \sup_{B(z,|x-y|)} u - \inf_{B(z,|x-y|)} u \,\le\, \theta \big( \frac{|x-y|}{r} \big) \inf_{B(z,|x-y|)} u \,\le \, \theta \big( \frac{|x-y|}{r} \big) \sup_{B_R}\, u \,.
$$
\qed

Even though this corollary gives merely the continuity estimates, 
%but the supremum is taken over the ball $B(x_0, R)$ and not the whold space $\R^d$ as in the existing literature
notice that the supremum is taken over the ball $B(x_0, R)$ and not the whole space $\R^d$ as in the existing literature
(see \cite{ BK, BKKu, B2, CSi, F, HK, Ka, Si, Sz2}).

\section{Relative Fatou  Theorem }

   In this section, we assume that $d \ge 2$. In the case $d= 2$,
 we will always assume the following:
 
 \medskip
\noindent
{\bf (A2)} : {\it
There
exists $\gamma\in (0, 1)$ such that
$
\liminf_{\lambda \to 0} \, {\phi(\lambda)}/{\lambda^{\gamma}}>0.
$
}
 \medskip

Then by the criterion of Chung-Fuchs type, the process $X$ is transient under this assumption (see \cite[(13.3.1)]{KSV2}).

In this section, using Theorem \ref{th:osc} we prove the relative Fatou theorem. 
The proofs of the results in this section are
similar to the corresponding parts of
\cite{K2}. For this reason, some proofs in this section will be omitted.

In this section, we assume that $D$ is a bounded $\kappa$-fat open set. We recall the definition of $\kappa$-fat open set.
   \begin{defn}\label{fat}
Let $\kappa \in (0,1/2]$. We say that an open set $D$ in $\R^d$ is
$\kappa$-fat if there exists $R>0$ such that for each $Q \in
\partial D$ and $r \in (0, R)$, $D \cap B(Q,r)$ contains a ball
$B(A_r(Q),\kappa r)$. The pair $(R, \kappa)$ is called the
characteristics of the $\kappa$-fat open set $D$.
\end{defn}
Note that all Lipschitz domains and all non-tangentially accessible
domains (see \cite{JK} for the definition) are $\kappa$-fat. The
boundary of a $\kappa$-fat open set may be not rectifiable, and in general, no regularity of its boundary can be inferred.
A bounded $\kappa$-fat open set may be disconnected.

The following boundary Harnack principle is the main result in \cite{ KSV, KSV2}.
\begin{thm}\emph{(\cite[Theorem 4.8]{KSV}, \cite[Theorem 13.4.22]{KSV2})}\label{BHP}
Suppose that $D$ is  a $\kappa$-fat open set with the
characteristics $(R, \kappa)$. There exists a  constant $c= c(\alpha, d, \ell, R,\kappa)>1$ such that if $r \le
R\wedge \frac14$ and $Q\in\partial D$, then for any nonnegative
functions $u, v$ in $\R^d$ which are regular harmonic in $D\cap B(Q,
2r)$ with respect to $X$ and vanish in $D^c \cap B(Q, 2r)$, we have
$$
c^{-1}\,\frac{u(A_r(Q))}{v(A_r(Q))}\,\le\, \frac{u(x)}{v(x)}\,\le
c \,\frac{u(A_r(Q))}{v(A_r(Q))} \qquad \text{for }\,x\in D\cap B(Q, \frac{r}2)\, .
$$
\end{thm}

 Let $x_0\in D$ be fixed and set
\[
M_D(x, y):=\frac{G_D(x, y)}{G_D(x_0, y)}, \qquad \text{for }\,x, y\in D~ \text{ and }\,
y\neq x_0.
\]
For each fixed $z\in \partial D$ and $x\in D$,
let
$
M_D(x, z):=\lim_{D\ni y\to z}M_D(x,y),
$
which exists by \cite[Theorem 5.5]{KSV}. For each $z\in
\partial D$, set $M_D(x, z)$
to be zero for $x\in D^c$. $M_D$ is called the Martin kernel of $D$ with respect to $X$.

As a consequence of  \cite[Theorem 5.11]{KSV},
for every nonnegative harmonic function $h$ for $X^D$, there exists a unique
finite measure $\nu$ on $\partial D$ such that
$$
h(x)=\int_{\partial D}M_D(x, z)\nu(dz) \qquad \text{for }\,x\in D.
$$
$\nu$ is called the Martin measure of $h$.

We will use $G(x, y)=G(x-y) = \int_0 ^{\ift} p(t,x,y)\,dt$
to denote the Green function of $X$.
$G$ is
radially decreasing and continuous
in $\R^d\setminus \{0\}$.

\medskip
The proof of the next result is similar to \cite[Theorem 2.4]{CS3} and \cite[Lemma 3.2]{K2}.
\begin{lemma}\label{lemma:boy}
For each $z \in \partial D$,
 $M_D(\, \cdot \, , z)$ is bounded regular harmonic
in $D \setminus B(z, \eps)$ for every $\eps > 0$.
\end{lemma}

\pf
Fix $z \in \partial D$ and $\eps > 0$, and let $h(x):=M_D(x, z)$ for $x \in \R^d$.
Note that $G(x,y) \ge G_D (x,y) $. By \cite[Theorem 13.3.2]{KSV2}, \cite[Lemma 3.3]{KSV3} and Theorem \ref{BHP}, there exist $c_1, c_2 >0$ which depend on  $ \alpha, d, \ell, \kp, R$ and  $ \text{diam} (D)$ such that for every $x \in D \setminus B(z, \eps/2)$,
\begin{eqnarray*}
&&h(x) = M_D(x,z)\,=\,\lim_{D \ni y \to z} \frac{G_D(x,y)} {G_D(x_0,y)}\,
\le\, c_1 \,  \frac{G_D(x,A)} {G_D(x_0,A)} \nonumber \\
&&\le\, c_1  \, \frac{G(x,A)} {G_D(x_0,A)}
\,\le\, c_2 \,\sup_{y \in D \setminus B(z, \eps/2)}  \frac1{|y-A|^{d} \, \phi(|y-A|^{-2})G_D(x_0,A)} < \infty 
\end{eqnarray*}
where $A:=A_{\eps/16}(z)$ (see Definition \ref{fat}).
Take an increasing sequence of smooth open sets $\{ D_m \}_{m \ge 1}$ such that $\overline{D_m} \subset D_{m+1}$ and $\cup^{\infty}_{m=1} D_m = D \setminus B(z, \eps)$.
Set $\tau_m := \tau_{D_m}$ and $\tau_{\infty} := \tau_{D \setminus B(z, \eps)}$ .
Then $\tau_m \uparrow \tau_{\infty}$ and
$\lim_{m \to \infty} X_{\tau_m} = X_{\tau_\infty} $ by quasi-left continuity of $X$.
Set $E = \{ ~\tau_m = \tau_\infty  ~\mbox{ for some } m \ge 1 \}$ and $N$ be the set of irregular boundary points of $D$. Since $X$ is symmetric, by \cite[(VI.4.6), (VI.4.10)]{BG} we get
\begin{equation}\label{zero1}
\P_x(X_{\tau_\infty} \in N)=0  ~~~ \text{for }\, x \in D.
\end{equation}
We also know from \cite[Lemma 5.9(i)]{KSV} that if $w \in \partial D, w \not= z $
and $w$ is a regular boundary point, then
$
h(x) \to 0 $ as  $ x \to w$
so that $h$ is continuous on $\overline{D \setminus B(z, \eps)} \setminus N$.
Since $h$ is bounded on $\R^d \setminus B(z, \eps/2)$,
by the bounded convergence theorem and (\ref{zero1}), we have
\begin{eqnarray} \label{e:m_bdd}
&&\lim_{m \to \infty} \E_x \left[ \,h (X_{\tau_m}) \,;\, \tau_m < \tau_\infty  \right]
=\lim_{m \to \infty} \E_x \left[ \,h (X_{\tau_m})1_{ \overline{D \setminus B(z, \eps)} \setminus N}(X_{\tau_m}) \,;\, \tau_m < \tau_\infty \right] \nonumber \\
&&= \E_x \left[ \,h (X_{\tau_\infty})1_{\overline{D \setminus B(z, \eps)} \setminus N}(X_{\tau_\infty}) \,;\, E\,^c\,\right]= \E_x \left[ \,h (X_{\tau_\infty}) \,;\, E \,^c\,\right].
\end{eqnarray}
Since $\tau_m \uparrow \tau_\ift$ and $\{\tau_m = \tau_\ift \}= \{\tau_n = \tau_\ift,\, n\ge m \}\uparrow
E$ as $m\rightarrow \ift$, by \eqref{e:m_bdd} and the monotone convergence theorem,
\begin{eqnarray*}
h(x) \,=\, \lim_{m \rightarrow \ift} \E_x [h(X_{\tau_m})]  &=& \lim_{m \rightarrow \ift} \E_x [h(X_{\tau_m})\,;\, \tau_m < \tau_\infty] +\lim_{m \rightarrow \ift} \E_x [h(X_{\tau_\ift})\,;\, \tau_m = \tau_\infty] \\
&=&  \E_x [h(X_{\tau_{\ift}})\, ; \,E\,^c \,] +  \E_x [h(X_{\tau_\ift}) \,;\, E \,] \,=\, \E_x [h(X_{\tau_{\ift}})]\,.
\end{eqnarray*}\qed

Throughout this paper, ${\cal F}_t$ is augmented right continuous $\sigma$-fields generated by $X_t ^D$. For a positive harmonic  function $h$ with respect to $X^D$, we let $(\P^h_x , X^h_t)$ be the $h$-transform of  $(\P_x , X^D_t)$, that is,
$$\P^h_x (A) := \E_x \left[ \frac{h(X^D_t)}{h(x)}; A \right]\quad \mbox{ if } A \in {\cal F}_t \, . $$
When $h(\cdot)=M_D (\cdot, z)$, we use the notation $(\P^z_x , X^z_t) := (\P^h_x , X^h_t)$ so that
 $(\P^z_x , X^z_t)$ is  $M_D(\cdot , z)$-transform of $(\P_x , X^D_t)$.

Let $\tau_D ^z$ be the life time of $X^z$. Using \cite[Theorem 3.10]{KPS} and {\bf (A1)}, the proof of the next result is similar to \cite[Theorem 3.3]{K2}.
\begin{thm}\label{thm:son1}
$$
\P^z_x \Big( \lim_{t \uparrow \tau^z_D} X^z_t = z ,~ \tau^z_D< \infty \Big) = 1 ~~~\mbox{ for every } x \in D,~ z \in \partial D.
$$
\end{thm}

\pf
See \cite[Theorem 3.3]{K2}.
\qed

The  following result is a simple consequence of Theorem \ref{thm:son1}.

\begin{prop}\label{prop:3.3}
Let $h$ be a positive harmonic function with respect to $X^D$ with Martin measure $\nu$ .
Then
\begin{equation*}
\P^h_{x}\Big(A \cap \big\{\lim_{t \uparrow \tau^h_D} X^h_t \in K\big\} \Big)=\frac1{h(x)} \int_K  M_D(x,z)\P^z_{x} (A)\nu(dz)
\end{equation*}
for every $x \in D$, $A \in {\cal F}_{\tau_D}$ and  Borel subset $K$ of $\partial D$.
\end{prop}
\pf
See \cite[Proposition 3.5]{K2}.
\qed

\begin{defn}\label{def:3.4}
$A \in {\cal F}_{\tau_D}$ is shift-invariant if whenever $T < \tau_D $ is a stopping time, $ 1_A \circ \theta_T = 1_A$ $\P_x$-a.s. for every $x \in D$.
\end{defn}

Using
\cite[Theorem 5.11]{KSV}, the proof of the next proposition is the same as the one in \cite[Proposition 3.7]{K2} (see also \cite[page 196]{B}).

\begin{prop}\label{prop:3.5}
\emph{(0-1 law)} If $A$ is shift-invariant, then $x \to \P^z_x (A)$ is a constant function which is either $0$ or $1$.
\end{prop}

Using \eqref{e:asumptiononle}, \cite[Theorem 1.5.3]{BGT} and the 0-version of \cite[Theorem 1.5.11]{BGT}, we have  the following inequalities; there exists  $c=c(\af, d, \ell)
>0$
such that
\bee
s^{d}{\phi(s^{-2})} \,\le \, c \,
{r^{d}}{\phi(r^{-2})} \qquad \text{for  }\,
0<s<r\le 4 \\ \label{el9}
\eee
and
\bee
 \quad \int_{0} ^r \frac{1}{s\, \phi(s^{-2})}ds
 \,\le \, c \,
\frac{1}{\phi(r^{-2})} \qquad \text{for  }\,  0<r\le 4. \label{el8}
\eee

From now on, we use notations $T_B : = \inf \{ t>0 :\, X_t \in B \}$, $T^z _B := \inf\{t>0: \, X^z_t \in B \}$ and $B^{\lambda}_y := B(y, \lambda \delta_D (y))$ for the convenience.
\begin{prop}\label{prop:3.2}
There exists $c = c (  \alpha, \ell, D)>1$ such that if $ 0 < \lambda < 1/2$ and  $x, y \in D $ with $|y - x | > 2 \delta_D (y) $,  then
\begin{equation*}
 \P_{x} \left(T_{B^{\lambda}_y}  < \tau_D \right) \ge c ~G_D(x, y) \lambda^{d} \delta_D (y)^{d} \phi\big( (2 \lambda \delta_D (y))^{-2} \big) \, .
\end{equation*}
\end{prop}
\pf
Fix $\ld \in (0,1/2)$ and $x, y \in D $ with $|y - x | > 2 \delta_D (y) $.
Since $x \not\in B(y, \delta_D(y))$, by \cite[Theorem 2.14]{KSV3} we get
\begin{equation}\label{eqn:3.6}
\E_x \Big[\int^{\tau_D}_0 1_{B^{\lambda}_y} (X_s) ds \Big]= \int_{B^{\lambda}_y}G_D (x,z)dz \,\ge\, c_1 G_D(x , y) \lambda^d\delta_D (y)^d \,.
\end{equation}
On the other hand, by the strong Markov property,
\begin{eqnarray}\label{eqn:3.7}
\E_x \Big[\int^{\tau_D}_0 1_{B^{\lambda}_y} (X_s) ds \Big]&=&
\E_x \left[\E_{X_{T_{B^{\lambda}_y}}} \Big[\int^{\tau_D}_0 1_{B^{\lambda}_y} (X_s) ds \Big] : T_{B^{\lambda}_y} < \tau_D \right]\nonumber \\
&\le& \P_{x}\left(T_{B^{\lambda}_y} <\tau_D \right) \sup_{w \in \overline{B^{\lambda}_y}} \E_w \Big[\int^{\tau_D}_0 1_{B^{\lambda}_y} (X_s) ds \Big]\, .
\end{eqnarray}
Note that since $ 0 < \lambda \, \delta_D (y) \le \text{diam}(D)$,
by \eqref{el8} and \cite[Theorem 13.3.2]{KSV2}, we obtain for every $w \in \overline{{B^{\lambda}_y}}$
\begin{eqnarray*}
&&\E_w \Big[\int^{\tau_D}_0 1_{B^{\lambda}_y} (X_s) ds \Big]\,\le\, \int_{B^{\lambda}_y} G(w-v)dv \,\le\, c_2 \int_{B^{\lambda}_y} \frac{dv}{|w-v|^{d}\phi(|w-v|^{-2})}\\
&&\le c_2
\int_{\{|w-v|\le 2 \lambda \delta_D (y) \}} \frac{dv}{|w-v|^{d}\phi(|w-v|^{-2})}\,=\,c_3 \int_{0}^{ 2 \lambda \delta_D (y)} \frac{1}{s \phi(s^{-2})} ds
 \,\le \, c_4 \,
\frac{1}{\phi\big( (2 \lambda \delta_D (y))^{-2}\big)} \, .
\end{eqnarray*}
Combining this with \eqref{eqn:3.6}--\eqref{eqn:3.7}, we finish the proof. 
\qed

Now we define the Stolz open set for $\kappa$-fat open set $D$ with the characteristics $(R,\kappa)$.
\begin{defn}\label{Stolz}
For $z \in \partial D$ and $\beta > (1-\kappa)/\kappa$, let
$
A^{\beta}_z := \{ y \in D \,;\, \delta_D(y) < R \,\wedge\, \big( \delta_D (x_0) /3 \big) \mbox{ and } |y-z| <  \beta\, \delta_D(y) \}.
$
We call $A^{\beta}_z$ the Stolz open set for $D$ at $z$ with the angle $\beta $.
\end{defn}

Since $\beta > (1-\kappa)/\kappa$,  there exists a sequence $\{y_k\}_{k \ge 1} \subset A^{\beta}_z$
such that $\lim_{k \to \infty} y_k =z$ (see \cite[Lemma 3.9]{K2}).

\begin{prop}\label{prop:3.6}
Given $\beta > (1-\kappa)/\kappa$ and $x\in D $, there exists $c = c( \alpha , \beta, D, x) > 0$ such that for every $z \in \partial D$, $\lambda \in (0, 1/2)$  and $y \in A^{\beta}_z$ with
$\delta_D (y) \le \frac12|x-y| \wedge  \delta_D (x)$,
 we have
$$
\P^z_{x} \Big(T^z_{B^{\lambda}_y} < \tau^z_D \Big) ~>~ c \,\lambda^{d}\, \frac{  \phi \big( (2 \lambda \delta_D (y))^{-2} \big)
}      {   \phi \big( (\delta_D (y)/8)^{-2} \big)} .
$$
\end{prop}
\pf
Fix $\beta > (1-\kappa)/\kappa$, $z \in \partial D$, $x\in D $, $\lambda \in (0, 1/2)$ and $y \in A^{\beta}_z$ with
$\delta_D (y) \le \frac12|x-y| \wedge  \delta_D (x)$. Let $z_1:=A_{\delta_D (y)/8}(z)$ so that
$
B(z_1, \kappa \, \delta_D (y) /8)
\subset
B(z, \delta_D (y)/8) \cap D
$
 and  fix $z_2 \in \partial B(y, \delta_D (y)/8)$.
Since $M_D(\cdot,z)$ is a harmonic function with respect to $X$ in $D$ (Lemma \ref{lemma:boy}), by Harnack principle (\cite[Theorem 2.14]{KSV3}) and Proposition \ref{prop:3.2} we have
\begin{eqnarray*}
\P^z_{x} \left( T^z_{B^{\lambda}_y} < \tau^z_D \right)
&=& \E_{x} \left[ \frac{M_D(X_{T_{B^{\lambda}_y}} , z)}{M_D (x,z)} \,;\,  T_{B^{\lambda}_y} < \tau_D \right]\,\geq\, c_1 ~ \P_{x} \left(  T_{B^{\lambda}_y} < \tau_D \right)~\frac{ M_D(y,z)}{ M_D(x,z)}\\
&\geq& c_2 \, G_D(x, y) \lambda^{d} \delta_D (y)^{d}\phi \big((2 \lambda \delta_D (y))^{-2} \big)
 \lim_{D \ni w  \rightarrow z} \frac{G_D(y,w)}{G_D(x,w)} \\
& \ge& c_3 \, G_D(x, y) \lambda^{d} \delta_D (y)^{d}\phi \big((2 \lambda \delta_D (y))^{-2}\big)
 \frac{G_D(y,z_1)}{G_D(x,z_1)} .
\end{eqnarray*}
The last inequality comes from Theorem \ref{BHP} because $|y-z| \wedge |x-z| \,>\, \delta_D (y)/2$.
We see that
$ \delta_D(z_1) \ge {\kappa \delta_D (y)}/{8} >{\delta_D (y)}/{\big(8(\beta+1)\big)}$,
$\delta_D(z_2) >  {\delta_D (y)}/{2}$ and $
|z_2- y| ={\delta_D (y)}/{8}.$ Moreover using our assumptions that
$\delta_D(y) \le  \delta_D(x)$ and $|x-y|\ge 2 \delta_D (y)$,
we have
\begin{eqnarray*}
|z_2-x| &\ge& |x-y|-|y-z_2| \,\ge\,  2\delta_D (y) -\frac{\delta_D (y)}{8}\,>\, \delta_D (y) \, ,\\
|z_1-x| &\ge& |x-z|-|z-z_1| \,\ge\, \delta_D (x) -\frac{\delta_D (y)}{8}
 \,>\, \frac{\delta_D (y)}{2}  
\end{eqnarray*}
and
\begin{equation*}  |z_1-y| \,\ge\,  |y-z|-|z_1-z| \,\ge\, \delta_D (y)- \frac{\delta_D (y)} {8} \,>\,  \frac{\delta_D (y)}{2} \, .
\end{equation*}
Thus $G_D(y,\cdot)$ and
$G_D(x,\cdot)$ are  harmonic functions  in
$
B(z_1,8^{-1}(\beta+1)^{-1}\delta_D (y))
\cup
B(z_2,8^{-1}(\beta+1)^{-1}\delta_D (y)).
$
Since
$
|z_1-z_2| \le |z_1-z|+|z-y|+|y-z_2|
< (4^{-1}+\beta) \, \delta_D (y)
$, by \cite[Theorem 2.14]{KSV3} we have
$G_D(y,z_1) \,\ge\, c_4 G_D(y,z_2)$ and
$G_D(x,z_1)\,\le\, c_5 G_D(x,z_2) \,\le\, c_6 G_D(x,y)$.
On the other hand,  by \cite[Lemma 3.3]{KSV3} and \eqref{el9},
we get $$
G_D(y, z_2) \,\ge\, c_7 \frac1{|y- z_2|^{d}\phi(|y-z_2|^{-2})} \,\ge\, c_8 \frac1{\delta_D (y)^{d}\phi\big( (\delta_D (y)/8)^{-2} \big)} .
$$
Combining these observations, we prove the proposition.
\qed

Now we are ready to show relative Fatou  theorem for harmonic function with respect to $X$ in $D$. The proof is similar to the proof of \cite[Theorem 3.13]{K2}. But, since we state a slightly more general version,  
we spell out detail for the reader's convenience.

\begin{thm}\label{T:Fatou}
Let $h$ be a positive harmonic  function with respect to $X^D$ with the Martin measure $\nu$.
If $u$ is a nonnegative function which is harmonic in $D$ with respect to $X$ and $x \in D$,  then for $\nu$-a.e. $z \in \partial D$, $\lim_{t \uparrow \tau^z_D} u(X^z_t)/h(X^z_t)$ exists and is finite $\P^z_{x}$-a.s.. Moreover, for every $x \in D$ and every $\bt> \frac{1-\kp}{\kp}$,
\begin{equation}\label{eqn:3.9n}
\lim_{t \uparrow \tau^z_D} \frac{u(X^z_t)}{h(X^z_t)}
=\lim_{ A^{\beta}_z \ni y \rightarrow z} \frac{u(y)}{h(y)} ~~~~\P^z_{x} \mbox{-a.s.}.
\end{equation}
In particular, for $\nu$-a.e. $z \in \partial D$,
\begin{equation}\label{eqn:3.9}
\lim_{ A^{\beta}_z \ni y \rightarrow z} \frac{u(y)}{h(y)} \mbox{ exists for every } \beta >\frac{1-\kappa}{\kappa}.
\end{equation}
\end{thm}

\pf
Without loss of generality, we assume  $\nu(\partial D) = 1$ and fix $x \in D$.
Note that $u$ is a non-negative and continuous superharmonic function with respect to $X^D$, i.e.,
for  $x\in B$, $
u(x)\ge \E_x\left[  u(X^D_{\tau_{B}})\right]
$
for every open set $B$ whose closure is a compact subset
of $D$.
Since $X^D$ is a Hunt process and $u$ is  non-negative and continuous superharmonic with respect to $X^D$, $u$ is excessive with respect to $X^D$
(see \cite[Corollary II.5.3]{BG} and the second part of the proof of \cite[Proposition II.6.7]{B}).
In particular, $\E_w[  u(X^D_t)] \le u(w)$ for every $w \in D$.
So by Markov property for conditional process (for example, see \cite[Chapter 11]{CW}), we have for every $t, s >0$
$$
\E^h_{x}\left[  \frac{u(X^h_{t+s})}{h(X^h_{t+s})} \,\big|\, {\cal F}_s\right]
= \E^h_{X^h_s}\left[  \frac{u(X^h_{t})}{h(X^h_{t})}\right]
=\frac1{h(X^h_s)} \, \E_{X^h_s}\left[  u(X^D_{t})\right]\le \frac{u(X^h_s)}{h(X^h_s)} \, .
$$
Therefore we see that $u(X^h_t)/h(X^h_t)$ is a non-negative supermartingale with respect to $\P^h_{x}$, and so the martingale convergence theorem gives
$\lim_{t \uparrow \tau^h_D} {u(X^h_t)}/{h(X^h_t)}$ exists and is finite $\P^h_{x}\mbox{-a.s.}$. Thus by Proposition \ref{prop:3.3}, for $\nu$-a.e. $z \in \partial D$,
\begin{equation}\label{eqn:3.10}
\P^z_{x} \left(\lim_{t \uparrow \tau^z_D} \frac{u(X^z_t)}{h(X^z_t)} \mbox { exists and is finite} \right) = 1.
\end{equation}
Fix $z \in \partial D$ satisfying \eqref{eqn:3.10} and  $\beta >(1-\kappa)/\kappa$.
By \eqref{e:asumptiononle} and Proposition \ref{prop:3.6}, for every  sequence $\{ y_k \}^{\infty}_{k=1} \subset A_z^{\beta} $ converging to $z$, $\ds \P^z_{x}\Big(T^z_{B^{\lambda}_{y_k}} < \tau^z_D ~\mbox { i.o.}\Big)\geq \liminf_{k \rightarrow \infty} \P^z_{x} \Big(T^z_{B^{\lambda}_{y_k}} < \tau^z_D\Big) >0 \quad \text{for every } \lambda \in (0, 1/2)  $. Since $\{ T^z_{B^{\lambda}_{y_k}} < \tau^z_D ~\mbox { i.o.} \}$ is shift-invariant, by Proposition \ref{prop:3.5},
\begin{equation}\label{eqn:3.12}
\P^z_{x} \left( X^z_t \mbox{ hits infinitely many } B^{\lambda}_{y_k} \right)
~=~\P^z_{x} \left(T^z_{B^{\lambda}_{y_k}} < \tau^z_D ~~\mbox { i.o.} \right)
~=~1  \quad \text{for every } \lambda \in (0, 1/2).
\end{equation}

Now let
$$m := \liminf_{ A^{\beta}_z \ni y \rightarrow z} \frac{u(y)}{h(y)}\quad  \text{ and }\quad  l := \limsup_{ A^{\beta}_z \ni y \rightarrow z} \frac{u(y)}{h(y)} .$$
First we note that
$l < \infty$. If not,
for any $M > 1$, there exists a sequence $\{ x_k \}^{\infty}_{k=1} \subset A_z^{\beta} $ such that $u(x_k)/h(x_k)> 4 M$ and $x_k \rightarrow z.$
By Theorem \ref{th:osc}, there exists $\lambda_1 = \lambda_1 (M, \alpha, d, \ell) >0 $ such that
${u(w)}/{h(w)} \geq {M^2(M+1)^{-2}u(x_k)}/{ h(x_k)} > M$
for every $ w \in B^{\lambda_1}_{x_k}$.
Thus by \eqref{eqn:3.12} we have
$
\lim_{t \uparrow \tau^z_D} {u(X^z_t)}/{h(X^z_t)} > M, \,\, \P^z_{x}$-a.s. for every $M >1$, which is a contradiction to \eqref{eqn:3.10}. Also if $l =0$, then $0 \le m \le l = 0$ so the theorem is clear. So we assume $0<l<\ift$.

For given $\eps > 0$, choose  sequences $\{ y_k \}^{\infty}_{k=1} \cup \{ z_k \}^{\infty}_{k=1}\subset A_z^{\beta} $ such that $u(y_k)/h(y_k) > (1+\eps )^{-1} l$, $ u(z_k)/h(z_k)< m+\eps $ and $y_k, z_k \rightarrow z.$
By Theorem \ref{th:osc}, there is $\lambda_2 = \lambda_2 (\eps, \alpha, d,\ell)>0 $ such that
\bee
\frac{u(w)}{h(w)} ~\geq~ \frac{u(y_k)}{(1+\eps )^2 \,h(y_k)}~ >~ \frac{l}{(1+\eps )^3} \quad \text{for every }w \in B^{\lambda_2}_{y_k}  \label{eqn:3.11} 
\eee
and 
\bee \frac{u(w)}{h(w)} ~\leq~ (1+\eps )^2\,\frac{u(z_k)} {h(z_k)} ~<~ (1+\eps )^2(m+\eps) \quad \text{for every } w \in B^{\lambda_2}_{z_k}. \label{eqn:3.11A}
\eee
Applying \eqref{eqn:3.10}--\eqref{eqn:3.12} to \eqref{eqn:3.11}--\eqref{eqn:3.11A} and letting $\eps \downarrow 0$, we obtain both \eqref{eqn:3.9n} and \eqref{eqn:3.9}.
\qed

If $u$ and $h$ are harmonic  functions in $D$ and $u/h$ is bounded, then $u$ can be recovered from non-tangential boundary limit values of $u/h$.
\begin{thm}\label{rep}
If $u$ is a harmonic  function in $D$ with respect to $X$ and $u/h$ is bounded for a positive harmonic  function $h$ in $D$ with respect to $X^D$ with the Martin measure $\nu$, then for every $x \in D$
$$
u(x)~=~h(x)\, \E^h_x\Big[\varphi_u \Big(\lim_{t \uparrow \tau^h_D} X^h_t \Big)\Big]
$$
where $
\varphi_u(z):= \lim_{ A^{\beta}_z \ni x \rightarrow z} {u(x)}/{h(x)}$, $\beta >{(1-\kappa)}/{\kappa}$,
 which is well-defined for $\nu$-a.e. $z \in \partial D$. If we further assume that $u$ is   positive in $D$, then
 $\varphi_u(z)$ is
Radon-Nikodym derivative of the (unique) Martin measure $\mu_u$ with respect to $\nu$.
\end{thm}
\pf
Using our Propositions \ref{prop:3.3} and \ref{prop:3.5}, the proof is the same as \cite[Theorem 3.18]{K2} (There are typos in the proof of \cite[Theorem 3.18]{K2} ; $v$ should be replaced by $h$).
\qed

When the boundary of $D$ is sufficiently smooth, by \cite[Theorem 1.1]{KSV3} Martin kernel enjoys the following estimate:
\begin{equation}\label{ennnw1} 
c^{-1} \big(\phi  (\delta_D(x)^{-2})\big)^{-1/2}      \,|x-z|^{-d} \,\le\, M_D(x,z)
\le c\big(\phi  (\delta_D(x)^{-2})\big)^{-1/2}  \, |x-z|^{-d}\, .
\end{equation}

Now suppose that $d=2$, $D=B:=B(0,1)$, $x_0=0$ and $\sigma_1$ is the normalized surface measure on $\partial B$. It is showed in \cite{K2} that the Stolz domain is the best possible one for Fatou theorem in $B$ for $(-\Delta)^{\alpha/2}$-harmonic function. Similarly, using \eqref{ennnw1}, we can show that our Stolz open set is also the best possible one here.

A curve $\sC_0$ is called a tangential curve in $B$ which ends on $\partial B$ if $\sC_0 \cap \partial B = \{w_0\} \in \partial B$, $\sC_0 \setminus \{w_0\} \subset B$ and there are no $r > 0$ and $ \beta >1$ such that $\sC_0 \cap B(w_0,r) \subset A^{\beta}_{w_0} \cap B(w_0,r)$.

\begin{thm}\label{counter_thm}
Let
$ \ds h(x):=\int_{ \partial B } M_B(x, w) \sigma_1(dw)$, $\sC_0$ be a tangential curve in $B$ which ends on $\partial B$ and $\sC_{\theta}$ be the rotation of $\sC_{0}$ about $x_0$ through an angle $\theta.$
Then there exists a positive harmonic function $u$ with respect to $X$  in $B:=B(x_0,1)$ such that for a.e. $\theta \in [0, 2\pi]$ with respect to Lebesgue measure,
$$
\lim_{|x| \rightarrow 1, \, x \in \sC_{\theta}} \frac{u(x)}{h(x)} \mbox{ does not exist}.
$$
\end{thm}
\pf
See \cite[Lemma 3.22 and Theorem 3.23]{K2}.
\qed

With the relative Fatou theorem
given in Theorem \ref{T:Fatou}, the proof of Theorem \ref{T:Fatou_FK}
almost identical to the corresponding parts of
\cite{K2}. For this reason, the proof of Theorem \ref{T:Fatou_FK} will be omitted. We refer \cite{C, CK2, K2} for the definitions of $\S_\infty(X^D)$ and $\A_\infty (X^D)$

For a smooth measure $\mu$ associated with a continuous additive
functional $A^\mu$ and a Borel measurable function $F$ on $D\times D$ that
vanishes along the diagonal, define
$$
e_{A^\mu+ F}(t):=\exp \Big( A^\mu_t + \sum_{0<s\le t}F(X^D_{s-}, X^D_s)
\Big) \quad\text{for }\, t\geq 0.
$$
Let $\mu\in \S_\infty(X^D)$ and $F\in \A_\infty (X^D)$
such that the gauge function $x\mapsto \E_x \left[
e_{A^\mu+F}(\tau_D)\right]$ is bounded.
A  Borel measurable function
$k$ defined on $D$ is said to be a positive $(\mu, F)$-harmonic function if $k >0$ and
$ \E_x \big[ e_{A^\mu+F}(\tau_B) k(X^D_{\tau_B})
\big] \,=\, k(x)$
for every open set $B$ whose closure is a compact subset of $D$ and $x\in B$.
By \cite[Theorem 5.16 and Section 6]{CK2},
there is a unique finite measure $\nu$ on $\partial D$ such that
$
 k(x) =\int_{\partial D} K_D(x, z) \,\nu (dz)
$, where $K_D (x,z)$ is the Martin kernel for the semigroup $Q_t f(x) := \E_x [e_{A^\mu+F} (t) f(X_t^D)]$. We call $\nu$ the Martin-representing measure of $k$.

\begin{thm}\label{T:Fatou_FK}
Let $D$ be a bounded $\kappa$-fat open set and $k$ be a positive $(\mu , F)$-harmonic function with the Martin-representing measure $\nu$.
If $u$ is a nonnegative $(\mu, F)$-harmonic function, then for $\nu$-a.e. $z \in \partial D$,
$
\lim_{ A^{\beta}_z \ni x \rightarrow z} \frac{u(x)}{k(x)}
$ exists for every  $\beta >(1-\kappa)/{\kappa}.
$
\end{thm}
\pf
See the proof of \cite[Theorem 4.7]{K2}.
\qed

Using the same argument as the one in \cite[Lemma 4.9 and Theorem 4.10]{K2}, one can see that the Stolz open set is the best possible one like Theorem \ref{counter_thm}.

\vspace{.1in}
\begin{singlespace}

\end{singlespace}

{\bf Panki Kim}

Department of Mathematical Sciences and Research Institute of Mathematics,
Seoul National University,
Building 27, 1 Gwanak-ro, Gwanak-gu,
Seoul 151-747, Republic of Korea

E-mail: \texttt{pkim@snu.ac.kr}

\bigskip

{\bf Yunju Lee}

Department of Mathematical Sciences,
Seoul National University,
Building 27, 1 Gwanak-ro, Gwanak-gu,
Seoul 151-747, Republic of Korea

E-mail: \texttt{grape3@snu.ac.kr}

\bigskip

\end{document}